
\magnification=1200
\hfuzz=3pt
\overfullrule=0mm

\hsize=125mm
\hoffset=4mm


\font\tensymb=msam9
\font\fivesymb=msam5 at 5pt
\font\sevensymb=msam7  at 7pt
\newfam\symbfam
\scriptscriptfont\symbfam=\fivesymb
\textfont\symbfam=\tensymb
\scriptfont\symbfam=\sevensymb

\font\titlefont=cmbx10 at 15pt
\font\notefont=cmr8
\font\noteitfont=cmti8


\font\refttfont=cmtt10 at 9pt

\font\sc=cmcsc10 \rm


\def\Ker{{\rm Ker}}

\def\id{{\rm id}}

\def\ZZ{{\bf Z}}

\def\and{\quad\hbox{and}\quad}

\def\Pr{\noindent{\sc{Proof.--- }}}

\def\cqfd{ {\sevensymb {\char 3}}}

\def\hfl#1#2{\smash{\mathop{\hbox to 6mm{\rightarrowfill}}
\limits^{\scriptstyle#1}_{\scriptstyle#2}}}


\null
\noindent

\vskip 15pt

\centerline{\titlefont Homotopy formulas for}

\medskip

\centerline{\titlefont cyclic groups acting on rings}

\vskip 30pt

\centerline{\sc Christian Kassel}

\bigskip\bigskip
\noindent
{\sc Abstract.}
{\it The positive cohomology groups of a finite group acting on a ring 
vanish when the ring has a norm-one element. In this note we give
explicit homotopies on the level of cochains when the group is cyclic,
which allows us to express any cocycle of a cyclic group
as the coboundary of an explicit cochain. 
The formulas in this note are closely related to the effective problems considered in previous joint work 
with Eli Aljadeff.
}



\vskip 15pt
\bigskip

\noindent
Given a cyclic group acting on a (nonnecessarily commutative) ring,
we exhibit in this note homotopy formulas implying the vanishing of the cohomology 
of the group with coefficients in the ring when the latter has a norm-one
element. Using our formulas, we can express any cocycle of a cyclic
group as the coboundary of an explicit cochain.

Such formulas are closely related to the effective problems considered
in [1] and~[2]. Actually we discovered them when we tried to find an effective
version of [2, Proposition~3.1]. 
Our formulas may be added to the already long list\footnote{($^1$)}{\notefont
Mariusz Wodzicki wittily 
proposed the name {\noteitfont homologbook} for such a list}
of homotopy formulas that turn out to be useful 
in many areas of homological algebra, see~[4] for an extensive use 
of explicit homotopies in cyclic homology.
\medskip

Let $R$ be a ring (with unit~1),
$t$ be a ring automorphism of~$R$, 
and $n \geq 2$ be a natural number such that $t^n = \id_R$
(the identity of~$R$).
Define two $\ZZ$-linear endomorphisms $T, N: R \to R$  
by
$$T = t - \id_R \and N = \id_R + t + \cdots + t^{n-1}. $$
In the literature, the endomorphism $N$ is called the {\it norm map}
or the {\it trace map}. 
We have
$$t\circ N = N\circ t = N \and
T\circ N = N\circ T = 0. \eqno (1)$$
Let 
$R^t = \Ker(T) = \{ a\in R \; |\, t(a)  = a \}$
be the subgroup of $R$ consisting of all $t$-invariant elements of~$R$.
It is a subring of~$R$ and the endomorphisms $T$, $N$ are 
left and right $R^t$-linear.

\goodbreak
For any element $x$ of~$R$ define 
$\ZZ$-linear endomorphisms
$j_x$, $j'_x$, $h_x$, $h'_x$ of~$R$ by
$$j_x(a) = xa, \qquad j'_x(a) = xt(a),$$
$$h_x(a) = - \sum_{i= 1}^{n-1} \, t^i(x) \bigl( \id_R + t + \cdots + t^{i-1}\bigr)(a),$$
$$h'_x(a) = \sum_{i= 1}^{n-1} \, ( \id_R + t + \cdots + t^{i-1}) \bigr(xt^{-i}(a) \bigr),
$$
where $a\in R$. (By convention, $t^0 = \id_R$.)

\medskip\goodbreak
\noindent
{\bf Lemma~1}.---
{\it For all $x$, $a\in R$ we have
$$j'_x(a) - j_x(a) = xT(a) \and
h'_x(a) - h_x(a) = N(x) N(a) - N(xa).$$
}

\Pr
The first identity follows immediately from the definitions. 
For the second one we have
$$\eqalign{
h'_x(a) - h_x(a)
&= \sum_{k= 1}^{n-1} \sum_{i= 0}^{k-1}\, t^i(x)t^{i-k}(a) 
+ \sum_{i= 1}^{n-1} \sum_{j= 0}^{i-1}\, t^i(x)t^{j}(a) \cr
& = \sum_{0\leq i < k \leq n-1}\, t^i(x)t^{n+i-k}(a) 
+ \sum_{0\leq j < i \leq n-1} \, t^i(x)t^{j}(a). \cr
}$$
Set $j = n+i-k$ in the penultimate sum. 
Then $i < k \leq n-1$ is equivalent to $i < j \leq n-1$.
Therefore, 
$$\eqalign{
h'_x(a) - h_x(a)
& = \sum_{0\leq i < j \leq n-1}\, t^i(x)t^{j}(a) 
+ \sum_{0\leq j < i \leq n-1} \, t^i(x)t^{j}(a) \cr
& = \sum_{i, j \in \{0, \ldots, n-1\}\atop i\neq j}\, t^i(x)t^{j}(a) \cr
& = \sum_{i= 0}^{n-1} \sum_{j= 0}^{n-1}\, t^i(x)t^{j}(a) 
- \sum_{i= 0}^{n-1}\, t^i(x)t^{i}(a) \cr
& = \biggl(\sum_{i= 0}^{n-1}\, t^i(x) \biggr)
\biggl(\sum_{j= 0}^{n-1}t^{j}(a) \biggr)
- \sum_{i= 0}^{n-1}\, t^i(xa) \cr
& = N(x) N(a) - N(xa).\cr
}$$
\hfill\cqfd

\medskip\goodbreak
\noindent
{\bf Corollary~1}.---
{\it For all $x\in R$ we have 
$$j_x \circ N = j'_x \circ N, \qquad T \circ h_x = T \circ h'_x,$$
$$N \circ j_x + h_x \circ T = N \circ j'_x + h'_x \circ T.$$
}

\Pr
By (1) and Lemma~1, for all $a\in R$ we have
$$(j'_x - j_x)\bigl( N(a) \bigr) = x T\bigl( N(a) \bigr) = 0.$$
Similarly, 
$$T\bigl(h'_x(a) - h_x(a)\bigr) 
= t\bigl(N(x)\bigr) t\bigl(N(a)\bigr) - N(x) N(a) 
- T\bigl(N(xa)\bigr) = 0.$$
Finally,
$$(h'_x - h_x) \bigl(T(a)\bigr) =  N(x) N\bigl(T(a)\bigr) - N\bigl(xT(a)\bigr)
= - N \bigl( j'_x(a) - j_x(a) \bigr). \eqno \hbox{\cqfd}$$
\medskip

\goodbreak
We now state the main result of this note.

\medskip
\noindent
{\bf Proposition}.---
{\it For all $x\in R$ we have
$$N(x)\, \id_R = N \circ j_x + h_x \circ T, \eqno (2)$$
$$N(x)\, \id_R = j_x \circ N + T \circ h_x, \eqno (3)$$
$$N(x)\, \id_R = N \circ j'_x + h'_x \circ T, \eqno (4)$$
$$N(x)\, \id_R = j'_x \circ N + T \circ h'_x  \eqno (5)$$
where $N(x)\, \id_R$ is the left multiplication by $N(x)$ in~$R$.
}
\medskip

\Pr
By Corollary~1, Identity (4) (resp.~(3)) is equivalent to Identity~(2) (resp.~(5)).
It suffices to prove (2) and (5).

\smallskip
{\it Identity~(2):} We have
$$\eqalign{
h_x\bigl(T(a)\bigr) & = 
- \sum_{i= 1}^{n-1} \, t^i(x) 
\bigl( \id_R + t + \cdots + t^{i-1}\bigr) \bigl(t(a) - a \bigr) \cr
& = - \sum_{i= 1}^{n-1} \, t^i(x) \bigl(t^i(a) - a \bigr) \cr
& = - \sum_{i= 0}^{n-1} \, t^i(xa) + xa + \sum_{i=0}^{n-1} \, t^i(x) a - xa \cr
& =  - N(xa) + N(x) a\cr
& = N(x)a - N(j_x(a)).\cr
}$$

{\it Identity~(5):} We have
$$\eqalign{
T\bigl(h'_x(a)\bigr) 
& = \sum_{i= 1}^{n-1} \, (t-\id_R)
\Bigl( ( \id_R + t + \cdots + t^{i-1}) \bigr(xt^{-i}(a) \bigr) \Bigr) \cr
& = \sum_{i= 1}^{n-1} \, (t^i - \id_R) \bigl(xt^{-i}(a) \bigr) \cr
& = \sum_{i= 1}^{n-1} \, t^i\bigl(xt^{-i}(a)\bigr) - \sum_{i= 1}^{n-1} \, xt^{-i}(a)  \cr
& = \sum_{i= 0}^{n-1} \, t^i(x)a - xa -  \sum_{i=0}^{n-1} \, xt^{-i}(a) + xa \cr
& = N(x) a - xN(a) 
= N(x) a - x\, t \bigl(N(a)\bigr) \cr
& = N(x) a - j'_x\bigl(N(a)\bigr). \cr
}$$
\hfill\cqfd
\medskip

\goodbreak
The following {\it homotopy formulas} are immediate consequences of the proposition.

\medskip
\noindent
{\bf Corollary~2}.---
{\it If $x\in R$ satisfies $N(x) = 1$ (the unit of~$R$), then
$$N \circ j_x + h_x \circ T = N \circ j'_x + h'_x \circ T = \id_R,$$
$$j_x \circ N + T \circ h_x = j'_x \circ N + T \circ h'_x = \id_R. $$
}

Corollary~2 has the following interesting consequences. Suppose there is 
an element $x$ of~$R$ such that $N(x) = 1$ (this is equivalent to
the image of the norm map $N : R\to R$ being the subring~$R^t$).
Under this condition any element $a\in R^t$ (i.e., killed by~$T$) is the image
under~$N$ of explicit elements of~$R$, namely
$$a = N(xa) . \eqno (6)$$

Similarly, any element $a\in R$ killed by the norm map (i.~e., $N(a) = 0$)
is the image under~$T$ of explicit elements of~$R$, namely
$$a = T\bigl(h_x(a) \bigr) = T\bigl(h'_x(a) \bigr) . \eqno (7)$$
(The identity  $a = T\bigl(h'_x(a) \bigr)$
already appeared in [1, Lemma~1 and Formula~(3)].)

\medskip\goodbreak
\noindent
{\bf Concluding remarks}.
(i) Since $t^n = \id_R$, the automorphism $t$ induces a $G$-module structure on the ring~$R$,
where $G$ is the cyclic group of order~$n$.
Consider the (co)homology groups
$H^i(G,R)$ and $H_i(G,R)$ of $G$ with coefficients in this $G$-module.
As is well known (see, e.g.,  [3, Chap.~XII, \S~7]), 
these groups can be realized as the (co)homology groups of the periodic (co)chain complex
$$ \cdots \to R \, \hfl{T}{} \, R \, \hfl{N}{} \, R \, \hfl{T}{} \,
R \, \hfl{N}{} \, R \to \cdots . \eqno (8)$$

If there is an element $x\in R$ such that $N(x) = 1$, then 
by Corollary~2 the operators $j_x$, $j'_x$, $h_x$, $h'_x$
are homotopies for the complex~(8). 
Consequently, 
$$H^{i}(G,R) = 0 = H_{i}(G,R) \eqno (9)$$
for all $i>0$.
The vanishing of the (co)homology of a cyclic group with coefficients in a ring~$R$
in the presence of a norm-one element $x\in R$ has been observed in~[2]. 
(Actually, similar vanishing results hold for any finite group,
see [2, Proposition~3.1].)

With Formulas (6), (7) we can do better than~(9), namely 
we can express any (co)cycle in the complex~(8)
as the (co)boundary of an explicit (co)chain. 
Such effective formulas play an important role in [1] and ~[2, Sections~4--5].

(ii) If $R$ is uniquely $n$-divisible and $x = 1/n \in R$ is the unique element such that
$nx=1$, then $x$ is $t$-invariant and $N(x) = nx = 1$. Under this condition
the operators $h_x$ and $h'_x$ become (for $a\in R$)
$$h_x(a) = - {1 \over n}\, \sum_{j=0}^{n-1}\, (n-1-j)\, t^j(a)
\and
h'_x(a) =  {1 \over n}\, \sum_{j=0}^{n-1}\, j t^j(a). \eqno (10)$$

(iii) Let $R_n = \ZZ \langle t_0(X), t_1(X), t_2(X), \ldots, t_{n-1}(X) \rangle$ 
be the free ring on $n$ indeterminates $t_i(X)$ indexed by the cyclic group~$\ZZ/n$
of order~$n$. 
We consider the ring automorphism $t$ of $R_n$ determined by
$$t\bigl(t_i(X) \bigr) = t_{i+1}(X)$$
for all~$i\in \ZZ/n$. 
For any ring $R$ equipped with a ring automorphism~$t$ such that $t^n = \id_R$,
and any $x\in R$,
there is a unique ring map $f : R_n \to R$ 
that commutes with the automorphisms~$t$ and sends $X$ to~$x$. 
We may consider $R_n$ as the universal ring for the situation
considered in this note. The proposition holds in~$R_n$ (with $x$ replaced by~$X$). 

Corollary~2 holds in~$R_n$, but
also in the quotient of $R_n$ by the two-sided ideal
generated by $1 - \sum_{i\in \ZZ/n}\, t_i(X)$; this quotient-ring is universal
for all rings equipped with a $\ZZ/n$-action and an element of norm one.

In this sense the formulas in this note can be considered as universal.

\bigskip\bigskip\goodbreak
\centerline{\bf References}
\vskip 20pt

\noindent
[1] {\sc E.~Aljadeff, C.~Kassel},
{\it Explicit norm one elements for ring actions of finite abelian groups},
Israel J.~Math.\ 129 (2002), 99--108.
\smallskip

\noindent
[2] {\sc E.~Aljadeff, C.~Kassel},
{\it Norm formulas for finite groups and induction from elementary abelian subgroups},
31~pages,
arXiv:math.RA/0402145.
\smallskip

\noindent
[3] {\sc H.~Cartan, S.~Eilenberg},
{\it Homological algebra},
Princeton University Press, Princeton, 1956.
\smallskip

\noindent
[4] {\sc C.~Kassel},
{\it Homologie cyclique, caract\`ere de Chern et lemme de perturbation},
J.~reine angew.\ Math.\ 408 (1990), 159--180.
\smallskip

\vskip 30pt

\line{Christian Kassel \hfill}
\line{Institut de Recherche Math\'ematique Avanc\'ee \hfill}
\line{CNRS - Universit\'e  Louis Pasteur\hfill}
\line{7 rue Ren\'e Descartes \hfill}
\line{67084  Strasbourg Cedex, France \hfill}
\line{E-mail:  {\refttfont kassel@math.u-strasbg.fr}\hfill}

\bye